\documentclass[10pt]{article}
\usepackage{latexsym}
\usepackage{amsfonts}
\usepackage{amssymb}
\usepackage[active]{srcltx} 
\usepackage{color}

\newcommand\tstrut{\rule{0pt}{2.6ex}}

\newenvironment{Pf}{\noindent {\bf Proof.}}{\hfill$\Box$\vspace{5mm}}
\newcommand{\ba}{\begin{array}}\newcommand{\ea}{\end{array}}

\newcommand{\ns}{\rm}

\newcommand{\nse}{\kern-3pt\ns$=$}\newcommand{\qd}{\hfill$\Box$\medbreak}

\newcommand{\ext}{\raise1pt\hbox{$\ts\bigwedge$}}

\newcommand{\ts}{\textstyle}
\newcommand{\rf}[1]{(\ref{#1})}
\newcommand{\chii}{\raise2pt\hbox{$\chi$}}

\newcommand{\F}{\mbox{$\cal F$}}
\newcommand{\Fg}{\mbox{${\cal F}\kern-2pt_g$}}

\newcommand{\Mg}{\mbox{${\cal M}\kern-2pt_g$}}
\newcommand{\Ng}{\mbox{${\cal N}\kern-2pt_g$}}
\newcommand{\U}{\mbox{${\cal U}$}}\newcommand{\V}{V\kern-1pt}

\newcommand{\Gg}{\mbox{${\cal G}\kern-2pt_g$}}
\newcommand{\cir}{\raise1.6pt\hbox{\footnotesize$\circ$}}

\newcommand{\Res}[2]{\hbox{\ns Res}\kern-16pt\lower5pt\hbox{\footnotesize$_{#1}$}\kern2pt\left[#2\right]}
\newcommand{\qk}{quaternion-K\"ahler\kern2pt}\renewcommand{\,}{\kern1pt}

\newcommand{\Aut}{{\rm Aut}}
\newcommand{\aut}{{\rm aut}}

\newcommand{\End}{{\rm End}}

\newcommand{\dirac}{/\kern-5pt\partial}

\newcommand{\Spin}{Spin\,\,\,}

\newcommand{\Lie}{{\cal L}}

\newcommand{\lra}{\longrightarrow}
\renewcommand{\ts}{\textstyle}
0
\newtheorem{theo}{Theorem}[section]
\newtheorem{defi}{Definition}[section]

\newtheorem{prop}{Proposition}[section]\def\frac#1#2{{#1\over#2}}
\def\be#1\ee{\begin{equation}#1\end{equation}}

\def\be{\begin{equation}}
\def\ee{\end{equation}}

\def\Aut{\mathrm{Aut }}

\def\SU{\mathrm{SU}}
\def\U{\mathrm{U}}
\def\E{\mathrm{E}}
\def\F{\mathrm{F}}
\def\SO{\mathrm{SO}}

\def\End{\mathrm{End}}

\def\Sp{\mathrm{Sp}}
\def\Spin{\mathrm{Spin}}

\begin{document}
\title{Almost even-Clifford hermitian manifolds with large automorphism group}

\author{
Gerardo Arizmendi\footnote{Centro de
Investigaci\'on en Matem\'aticas, A. P. 402,
Guanajuato, Gto., C.P. 36000, M\'exico. E-mail: gerardo@cimat.mx}
\footnote{Partially supported by a
CONACYT scholarship},\,\,\,\,\,
Rafael Herrera\footnote{Centro de
Investigaci\'on en Matem\'aticas, A. P. 402,
Guanajuato, Gto., C.P. 36000, M\'exico. E-mail: rherrera@cimat.mx}
\footnote{Partially supported by
grants from CONACyT and LAISLA (CONACyT-CNRS)}
and Noemi Santana\footnote{Instituto de Matem\'aticas, UNAM, Unidad Cuernavaca, A.P.
6--60, C.P. 62131, Cuernavaca, Morelos, M\'exico. E-mail: noemi.santana@im.unam.mx}
\footnote{Partially
supported by 
grants of CONACyT and LAISLA (CONACyT-CNRS).}
}

\date{}

\maketitle

{
\abstract{

We study 
manifolds endowed with an (almost) even Clifford (hermitian) structure
and admitting a large automorphism group. We classify them
when they are simply connected and the dimension of
the automorphism group is maximal, and also prove a gap theorem for the dimension of the automorphism group.

}
}

\section{Introduction}

Recently, there has been some interest in manifolds admitting so-called even Clifford structures
\cite{Moroianu-Semmelmann,Moroianu-Pilca}. Here, we study such manifolds when they admit a large
automorphism group.
This type of problem has been studied on Riemannian manifolds
\cite{wang,wakakuwa}, almost
hermitian manifolds \cite{tanno}, and almost quaternion-hermitian manifolds \cite{Santana}.

It is a classical result \cite{kobayashi-1} that the maximal dimension of
the
isometry group of a connected $n$-dimensional Riemannian manifold is
${1\over 2}n(n+1)$. If the dimension is maximal, the manifold
is isometric to either
Euclidean space $\mathbb{R}^n$,
or the sphere $S^n$,
or projective space $\mathbb{RP}^n$,
or (simply connected) hyperbolic space.
Furthermore, in \cite{wang} it was shown that the
isometry group contains no $m$-dimensional closed subgroup where
\[{1\over 2} n(n-1)+1< m < {1\over 2} n(n+1). \]

In \cite{tanno}, it was shown that the automorphism group of a connected
$2n$-dimensional almost-hermitian manifold has dimension at most $n(n+2)$.
If the dimension of the automorphism is maximal, the manifold is isometric to either
complex Euclidean space $\mathbb{C}^n$, or
an open ball with K\"ahler structure with negative constant
holomorphic sectional curvature, or
complex projective space $\mathbb{CP}^n$.
\textcolor{black}{In this case, however, there is also a (unique) manifold whose automorphism group
has dimension one less than the maximal one, namely, euclidean space \cite{ishihara}.}

In \cite{Santana}, it was shown that the automorphism group of a connected
$4n$-dimensional almost quaternion-hermitian manifold has dimension at most $2n^2+5n+3$.
If the dimension of the automorphism group is between $2n^2+5n$ and $2n^2+5n+3$, the manifold
is isometric to either
quaternionic Euclidean space $\mathbb{H}^n$, or
quaternionic projective space $\mathbb{HP}^n$, or
quaternionic hyperbolic space.

In this paper, we will prove the analogous theorems for almost even-Clifford hermitian manifolds of rank
$r\geq 3$.
Our terminology differs from that of \cite{Moroianu-Semmelmann, Moroianu-Pilca} 
since we have added the words ``almost'' and ``hermitian'' since, in principle,  there
is no integrability condition on the structure and the compatibility with a Riemannian metric is an extra
condition. We
shall explore integrability conditions in the style of Gray \cite{Gray} in a
future paper. 

The note is organized as follows. In Section~\ref{sec: preliminaries} we
recall some preliminaries on Clifford algebras and representations, and almost even-Clifford hermitian
manifolds.
In Section~\ref{group}, we give and upper bound in the dimension of the automorphism group
(Proposition~\ref{dim-max}), determine the spaces whose automorphism group has dimension equal to this bound 
(Theorem~\ref{main-theorem}), and
prove a gap in the dimension of the automorphism group (Proposition~\ref{theo:gap theorem}).

{\bf Acknowledgments}. 
The second named author would like to thank  the 
International Centre for Theoretical Physics and the Institut des Hautes \'Etudes Scientifiques for their
hospitality and support.

\section{Preliminaries}\label{sec: preliminaries}

In this section we recall material that can also be consulted in \cite{Friedrich, Lawson}.
Let $Cl_r$ denote the Clifford algebra generated by all the products of the orthonormal vectors
$e_1, e_2, \ldots, e_r\in \mathbb{R}^r$ subject to the relations
\begin{eqnarray*}
e_j e_k + e_k e_j &=& -2\,\delta_{jk}, 
 \quad\mbox{for $1\leq j, k\leq r$.} 
\end{eqnarray*}
The even Clifford subalgebra $Cl_r^0$ is defined as the invariant (+1)-subspace of the involution of $Cl_r$
induced 
by the map $-{\rm Id}_{\mathbb{R}^r}$.
The Spin group $Spin(r)\subset Cl_r$ is the subset
\[Spin(r) =\{x_1x_2\cdots x_{2l-1}x_{2l}\,\,|\,\,x_j\in\mathbb{R}^r, \,\,
|x_j|=1,\,\,l\in\mathbb{N}\},\]
endowed with the product of the Clifford algebra.
The Lie algebra of $Spin(r)$ is
\[\mathfrak{spin}(r)=\mbox{span}\{e_ie_j\,\,|\,\,1\leq i< j \leq r\}.\]

Now, we summarize some results about real representations of $Cl_r^0$ in the next table (cf. \cite{Lawson}).
Here $d_r$ denotes the dimension of an irreducible representation of $Cl^0_r$ and $v_r$ the number of distinct
irreducible representations.
Let $\tilde\Delta_r$ denote the irreducible representation of $Cl_r^0$ for $r\not\equiv0$ $(\mbox{mod } 4) $
and $\tilde\Delta^{\pm}_r$ denote the irreducible representations for $r\equiv0$ $(\mbox{mod } 4)$.
\[\begin{array}{|c|c|c|c|c|c|}
\hline
r \mbox{ (mod 8)}&d_r&Cl_r^0&\tilde\Delta_r\,\,/\,\,\tilde\Delta_r^\pm\cong \mathbb{R}^{d_r}& v_r \tstrut\\
\hline
1&2^{\lfloor\frac r2\rfloor}&\mathbb R(d_r)&\mathbb{R}^{d_r}&1 \tstrut\\
\hline
2&2^{\frac r2}&\mathbb C(d_r/2)&\mathbb{C}^{d_r/2}&1 \tstrut\\
\hline
3&2^{\lfloor\frac r2\rfloor+1}&\mathbb H(d_r/4)&\mathbb{H}^{d_r/4}&1 \tstrut\\
\hline
4&2^{\frac r2}&\mathbb H(d_r/4)\oplus \mathbb H(d_r/4)&\mathbb{H}^{d_r/4}&2 \tstrut\\
\hline
5&2^{\lfloor\frac r2\rfloor+1}&\mathbb H(d_r/4)&\mathbb{H}^{d_r/4}&1 \tstrut\\
\hline
6&2^{\frac r2}&\mathbb C(d_r/2)&\mathbb{C}^{d_r/2}&1 \tstrut\\
\hline
7&2^{\lfloor\frac r2\rfloor}&\mathbb R(d_r)&\mathbb{R}^{d_r}&1 \tstrut\\
\hline
8&2^{\frac r2-1}&\mathbb R(d_r)\oplus \mathbb R(d_r)&\mathbb{R}^{d_r}&2 \tstrut\\
\hline
\end{array}
\]
\centerline{Table 1}

Note that
the representations are complex for $r\equiv 2,6$ $(\mbox{mod } 8)$ and quaternionic for $r\equiv 3,4,5$ 
$(\mbox{mod } 8)$.

\subsubsection*{Almost even-Clifford hermitian structures}

\begin{defi} 
 A {\em linear even-Clifford hermitian structure of rank $r$} on $\mathbb{R}^N$, $N\in\mathbb{N}$, is a
representation
\[Cl_r^0\longrightarrow \End(\mathbb{R}^N)\]
such that each bivector $e_ie_j$, $1\leq i< j \leq r$, is mapped to an antisymmetric endomorphism $J_{ij}$
satisfying
\begin{equation}
J_{ij}^2 = -{\rm Id}_{\mathbb{R}^N}.\label{eq:almost-complex-structures}
\end{equation}
\end{defi}

Notice that the subalgebra $\mathfrak{spin}(r)$ is mapped injectively into the skew-symmetric
endomorphisms $\End^-(\mathbb{R}^N)$. 

First, let us assume $r\not\equiv 0
\,\,\,({\rm mod}\,\,\,\, 4)$, $r>1$.
In this case, $\mathbb{R}^N$ decomposes into a sum of irreducible representations of $Cl_r^0$.
Since this algebra is simple, such irreducible representations can only be trivial or
copies of the standard representation
$\tilde\Delta_r$ of $Cl_r^0$ (cf. \cite{Lawson}). Due to
\rf{eq:almost-complex-structures}, there are no trivial summands in such a decomposition so that
\[\mathbb{R}^N =\tilde\Delta_r \otimes_{\mathbb{R}}\mathbb{R}^m\]
for some $m\in\mathbb{N}$.
Thus, we see that $\mathfrak{spin}(r)$ has an isomorphic image
\[\widehat{\mathfrak{spin}(r)}:=\mathfrak{spin}(r)\otimes \{{\rm Id}_{m\times m}\}\subset
\mathfrak{so}(d_rm).\]

Secondly, let us assume $r\equiv 0
\,\,\,({\rm mod}\,\,\,\, 4)$.
Recall that if $\hat\Delta_r$ is the irreducible representation of $Cl_r$, then by restricting this
representation to $Cl^0_r$ it splits as the sum of two inequivalent irreducible representations
\[\hat\Delta_r = \tilde{\Delta}_r^+ \oplus \tilde{\Delta}_r^-.\]
Since $\mathbb{R}^N$ is a representation of $Cl_r^0$ satisfying
\rf{eq:almost-complex-structures}, there are no trivial summands in such a decomposition so that
\begin{eqnarray*}
\mathbb{R}^N
&=& \tilde\Delta_r^+\otimes \mathbb{R}^{m_1} \oplus \tilde\Delta_r^-\otimes \mathbb{R}^{m_2},
\end{eqnarray*}
for some $m_1,m_2\in \mathbb{N}$.
By restricting this representation to $\mathfrak{spin}(r)\subset Cl_r^0$,
consider
\[\widehat{\mathfrak{spin}(r)}:=\mathfrak{spin}(r)^+\otimes \{{\rm Id}_{m_1\times m_1}\oplus
\mathbf{0}_{m_2\times m_2}\}\oplus
\mathfrak{spin}(r)^-\otimes\{\mathbf{0}_{m_1\times m_1}\oplus {\rm Id}_{m_2\times m_2}\}
\subset
\mathfrak{so}(d_rm_1+d_rm_2),\]
where $\mathfrak{spin}(r)^{\pm}$ are the images of $\mathfrak{spin}(r)$ in 
$\End(\tilde\Delta_r^{\pm})$
respectively.

\begin{defi} 
\begin{itemize}
 \item 
A {\em rank $r$ almost even-Clifford hermitian structure}, $r\geq 2$, on a Riemannian manifold
$M$
is a smoothly varying choice of linear even-Clifford hermitian structure on each tangent space of $M$.
Let $Q\subset \End^-(TM)$ denote the subbundle with fiber $\mathfrak{spin}(r)$.

\item A Riemannian manifold carrying such a structure will be called an {\em almost even-Clifford hermitian
manifold}.

\item  An almost even-Clifford hermitian structure on a Riemannian manifold $M$ is called {\em parallel} if 
the bundle $Q$ is  parallel with respect to the Levi-Civita connection on $M$.
\end{itemize}
\end{defi}

Notice that the definition of (parallel) even-Clifford structure in \cite{Moroianu-Semmelmann} implies the one
we have just given.

\section{Automorphism group}\label{group}

In this section we derive an upper bound on the dimension of the automorphism
group of an almost  even-Clifford hermitian manifold and classify the manifolds whose automorphism group's
dimension attains such an upper bound.

The {\em automorphism group}
of an almost even-Clifford hermitian manifold $M$, denoted by $\Aut(M)$,
is the (sub)group of isometries which preserve the almost even-Clifford hermitian structure.
A vector field $X$ on $M$ is an {\em infinitesimal automorphism} if it is a Killing
vector field that preserves the structure, i.e. locally
\[\mathcal{L}_X J_{ij} =\sum_{k<l} \alpha_{kl}^{(ij)}\, J_{kl},\]
for some (local) functions $\alpha_{kl}^{(ij)}$, 
where $\mathcal{L}_X$ denotes the Lie derivative in the direction of $X$. These vector fields form the Lie
algebra $\aut(M)$ of $\Aut(M)$.

\subsection{Upper bound}
Let $X$ be an infinitesimal automorphism of $M$. Consider
\begin{eqnarray*}
\Lie_X(J_{ij}(Y)) &=& (\Lie_XJ_{ij})(Y) + J_{ij}(\Lie_XY), 
\end{eqnarray*}
i.e.
\begin{eqnarray*}
\nabla_X(J_{ij}(Y))-\nabla_{J_{ij}(Y)}X
   &=& 
  \sum_{k<l} \alpha_{kl}^{(ij)}\, J_{kl}(Y)+J_{ij}(\nabla_XY-\nabla_YX).
\end{eqnarray*}
Now suppose we are calculations at a point $p$ where $X_p=0$, so that
\begin{eqnarray*}
-\nabla_{J_{ij}(Y)}X
   &=& 
  \sum_{k<l} \alpha_{kl}^{(ij)}\, J_{kl}(Y)-J_{ij}(\nabla_YX),
\end{eqnarray*}
i.e.
\begin{eqnarray*}
{}[J_{ij},\nabla X](Y)
   &=& 
  \sum_{k<l} \alpha_{kl}^{(ij)}\, J_{kl}(Y).
\end{eqnarray*}
Hence, $(\nabla X)_p$ is a skew-symmetric endomorphism such that
\begin{eqnarray*}
{}[J_{ij},\nabla X]
   &=& 
  \sum_{k<l} \alpha_{kl}^{(ij)}\, J_{kl}.
\end{eqnarray*}
i.e. $(\nabla X)_p$ belongs to the normalizer of
$\widehat{\mathfrak{spin}(r)}={\rm span}(J_{ij})$ in
$\End^-(T_pM)=\mathfrak{so}(T_pM)$. 
Such a normalizer has been calculated in \cite{Arizmendi-Herrera} and we list them for $r\geq 3$.
\[  
\begin{array}{|c|c|c|c|}
 \hline
  r \mbox{\ {\rm (mod 8)} }  &N
&N_{\mathfrak{so}(N)}(\widehat{\mathfrak{spin}({r})})&C_{\mathfrak{so}(N)}(\widehat{\mathfrak{spin}({r})})
\rule{0pt}{3ex}\\
\hline
 0 &
d_r(m_1+m_2)&\mathfrak{so}(m_1)\oplus\mathfrak{so}(m_2)\oplus\mathfrak{spin}(r)&\mathfrak{so}
(m_1)\oplus\mathfrak{so}(m_2)\tstrut\\
 \hline
 1 & d_rm&\mathfrak{so}(m)\oplus\mathfrak{spin}(r)&\mathfrak{so}(m)\tstrut\\
\hline
 2 & d_rm& \mathfrak{u}(m)\oplus\mathfrak{spin}(r)& \mathfrak{u}(m)\tstrut\\
\hline
 3 & d_rm& \mathfrak{sp}(m)\oplus\mathfrak{spin}(r)& \mathfrak{sp}(m)\tstrut\\
\hline
 4 & d_r(m_1+m_2)& \mathfrak{sp}(m_1)\oplus\mathfrak{sp}(m_2)\oplus\mathfrak{spin}(r)&
\mathfrak{sp}(m_1)\oplus\mathfrak{sp}(m_2)\tstrut\\
\hline
 5 & d_rm& \mathfrak{sp}(m)\oplus\mathfrak{spin}(r)& \mathfrak{sp}(m)\tstrut\\
\hline
 6 & d_rm& \mathfrak{u}(m)\oplus\mathfrak{spin}(r)& \mathfrak{u}(m)\tstrut\\
\hline
 7 & d_rm& \mathfrak{so}(m)\oplus\mathfrak{spin}(r)& \mathfrak{so}(m)\tstrut\\
\hline
  \end{array}
\]
\centerline{Table 2}

\begin{prop}\label{dim-max} Let $M$ be a $N$-dimensional almost even-Clifford hermitian
manifold.
Then
\[  
\begin{array}{|c|c|c|}
 \hline
  r \mbox{\ {\rm (mod 8)} }  &N=\dim(M) &\mbox{upper bound } d_{max}\geq \dim(\Aut(M)) 
\tstrut\\
\hline
 0 & d_r(m_1+m_2)& {m_1\choose 2}+{m_2\choose 2}+{r\choose
2}+d_r(m_1+m_2)\tstrut\\
 \hline
 1,7 & d_rm&{m\choose 2}+{r\choose 2}+d_rm\tstrut\\
\hline
 2,6 & d_rm& m^2+{r\choose 2}+d_rm\tstrut\\
\hline
 3,5 & d_rm& {2m+1\choose 2}+{r\choose 2}+d_rm\tstrut\\
\hline
 4 & d_r(m_1+m_2)& {2m_1+1\choose 2}+{2m_2+1\choose 2}+{r\choose
2}+d_r(m_1+m_2)\tstrut\\
\hline
  \end{array}
\]
\centerline{\rm Table 3}
\end{prop}
\qd
\subsection{Large automorphism group}

In this subsection, we determine the spaces that support an automorphism group of maximal dimension and prove
a gap in the dimension of the automorphism group.

\begin{prop}
Let $M$ be a $N$-dimensional, rank $r\geq 3$ almost even-Clifford
hermitian
manifold and assume that the dimension of its automorphism group is maximal. 
Then, for any $p\in M$, the
isotropy subgroup $A_p$ of $p$ is conjugate to
$C_{{SO}(N)}({{Spin}({r})})\cdot Spin(r)\subset SO(N)$.
\end{prop}
{\em Proof}.
The dimension of the orbit f $p$ under $\Aut(M)$ satisfies
\[\dim(\Aut(M))-\dim(A_p)\leq N,\]
so that
\begin{eqnarray*}
\dim(A_p)
  &\geq& \dim(\Aut(M))-N \\
  &=& d_{max}-N \\
  &=& \dim(C_{{SO}(N)}({{Spin}({r})})\cdot Spin(r)).
\end{eqnarray*}
The Lie algebra $\mathfrak{a}_p$ of $A_p$ maps one-to-one into
$C_{\mathfrak{so}(N)}(\widehat{\mathfrak{spin}({r})})\oplus \widehat{\mathfrak{spin}({r})}$ since a Killing
vector field $X$ is determined by its values $X_p$ and $(\nabla X)_p$. Hence, 
\[\mathfrak{a}_p \cong C_{\mathfrak{so}(N)}(\widehat{\mathfrak{spin}({r})})\oplus 
\widehat{\mathfrak{spin}({r})}.\]
\qd

\begin{prop}\label{prop:sym-parallel}
 Let $M$ be a rank $r\geq 3$ almost even-Clifford hermitian
manifold and assume that the dimension of its automorphism group is maximal. Then, $M$ is symmetric and the
almost even-Clifford hermitian structure is parallel.
\end{prop}
{\em Proof}.
Let $p\in M$ and $A_p$ denote its isotropy group. 
We know that $A_p=C_{SO(N)}(Spin(r))\cdot Spin(r)$.
Since $C_{SO(N)}(Spin(r))$ contains $1$ and $Spin(r)$ contains $-1$, we have $-1\in A_p$. 
Thus, there is an element $g\in A_p$ whose derivative
$dg_p=-{\rm Id}_{T_pM}$ in the isotropy representation of $A_p$ on $T_pM$. In other words, the automorphism
$g$ is a (global) symmetry at $p$ and $M$ is symmetric. 
Since these symmetries generate the translations along geodesics, $M$ has a
transitive group of automorphisms, not just isometries.

Proceeding as in \cite[p. 264]{Jost}, given a vector field $W$ with $W_p\not =0$, let $c(t)$ be the
geodesic with $\dot{c}(0) = W_p$ and $\tau_t$ be the group of translations along $c$.
Then \[Z_q:=\frac{d}{dt}\tau_t(q)_{|t=0}\] is an infitesimal automorphism since $\tau_t$ are automorphisms. 
We have that
$W_p=Z_p$.
For $v\in T_pM$, let $\gamma(s)$ be a curve with $\gamma'(0)=v$. Then
\begin{eqnarray*}
\nabla_vZ_p&=&\nabla_{\frac{\partial}{\partial s}}\frac{\partial}{\partial t}\tau_t(\gamma(s))_{|s=t=0}\\
&=&\nabla_{\frac{\partial}{\partial t}}\frac{\partial}{\partial s}\tau_t(\gamma(s))_{|s=t=0}\\
&=&\nabla_{\frac{\partial}{\partial t}}D\tau_t(v)_{|t=0}\\
&=&0,
\end{eqnarray*}
since  $D_{\tau_t}$ is parallel transport along $c$ and $D_{\tau_t}(v)$ is a parallel vector field along
$c$. Hence, for any vector $W_p\in T_pM$, we have an infinitesimal isometry $Z$ such that  
\[Z_p=W_p \quad\quad\mbox{and}\quad\quad  (\nabla Z)_p=0.\] 

Now, recall that 
\[\nabla_Z(Y) = \Lie_ZY + \nabla_YZ.\]
On the one hand,
\begin{eqnarray*}
(\nabla_W (J_{ij}(Y)))_p
   &=&
  (\nabla_ Z (J_{ij}(Y)))_p\\
   &=& 
  ((\nabla_ZJ_{ij})(Y))_p + (J_{ij}(\nabla_ZY))_p,
\end{eqnarray*}
and, on the other,
\begin{eqnarray*}
(\Lie_Z(J_{ij}(Y)) + \nabla_{J_{ij}(Y)}Z)_p
   &=&
((\Lie_ZJ_{ij})(Y))_p + (J_{ij}(\Lie_ZY))_p + (\nabla_{J_{ij}(Y)}Z)_p\\
   &=&
\left(\sum_{k<l}\alpha_{kl}^{(ij)}J_{kl}(Y)\right)_p + (J_{ij}(\nabla_ZY))_p- (J_{ij}(\nabla_YZ))_p +
(\nabla_{J_{ij}(Y)}Z)_p,
\end{eqnarray*}
so that
\begin{eqnarray*}
  ((\nabla_ZJ_{ij})(Y))_p 
   &=& 
\left(\sum_{k<l}\alpha_{kl}^{(ij)}J_{kl}(Y)\right)_p - (J_{ij}(\nabla_YZ))_p +
(\nabla_{J_{ij}(Y)}Z)_p\\
   &=&
\left(\sum_{k<l}\alpha_{kl}^{(ij)}J_{kl}(Y)\right)_p + [(\nabla Z)_p,(J_{ij})_p](Y_p)\\
   &=&
\left(\sum_{k<l}\alpha_{kl}^{(ij)}J_{kl}(Y)\right)_p ,
\end{eqnarray*}
i.e.
\begin{eqnarray*}
  (\nabla_WJ_{ij})_p 
   &=&
\left(\sum_{k<l}\alpha_{kl}^{(ij)}J_{kl}\right)_p .
\end{eqnarray*}
\qd

\begin{theo}\label{main-theorem}
 Let $M$ be a simply connected Riemannian almost even-Clifford hermitian manifold of rank $r\geq 3$ such that
the dimension of its group of automorphisms is maximal. Then $M$ is isometric to one of the
following spaces: 
\[\begin{array}{|c|c|}
\hline
r &M\\
\hline
\mbox{arbitrary}&\tilde\Delta_r^{\times m} \mbox{ or }(\tilde\Delta_r^+)^{\times
m_1}\oplus (\tilde\Delta_r^-)^{\times m_2}, \mbox{ for some $m,m_1,m_2\in\mathbb{N}$}
\tstrut\\
\hline
3& 
\Sp(k+1)/(\Sp(k)\times\Sp(1)),\Sp(k,1)/(\Sp(k)\times\Sp(1))
\tstrut\\
\hline
4&M_1\times M_2\mbox{, where }M_i
=(\tilde\Delta_3)^{\times m},\Sp(k+1)/(\Sp(k)\times\Sp(1)),\Sp(k,1)/(\Sp(k)\times\Sp(1))\tstrut\\
\hline
5&\Sp(k+2)/(\Sp(k)\times\Sp(2)),\Sp(k,2)/(\Sp(k)\times\Sp(2))\tstrut\\
\hline
6&\SU(k+4)/{\rm S}(\U(k)\times\U(4)),\SU(k,4)/{\rm S}(\U(k)\times\U(4))\tstrut\\
\hline
8&\SO(k+8)/(\SO(k)\times\SO(8)), SO(k,8)/(\SO(k)\times\SO(8))\tstrut\\
\hline
9&\F_4/\Spin(9),\F_4^{-20}/\Spin(9)\tstrut\\
\hline
10&\E_6/(\Spin(10)\cdot\U(1)), \E_6^{-14}/(\Spin(10)\cdot\U(1))\tstrut\\
\hline
12&\E_7/(\Spin(12)\cdot\SU(2)), \E_7^{-5}/(\Spin(12)\cdot\SU(2))\tstrut\\
\hline
16&\E_8/\Spin^+(16), \E_8^8/\Spin^+(16)\tstrut\\
\hline
\end{array}\]
\centerline{\rm Table 4}
\end{theo}
{\em Proof}. The flat case is clear and 
the case $r=3$ was dealt with in \cite{Santana}, 

For $r=4$, by \cite{Moroianu-Semmelmann}, $M$ is a Riemannian product $M_1\times M_2$ of quaternion-K\"ahler
manifolds.
We claim that
\[\dim(\Aut(M))
=\dim(\Aut(M_1))+\dim( \Aut(M_2)).\] 
Indeed, let $X\in \aut(M)$ an infinitesimal automorphism of $M$,
$X=X_1+X_2$ with $X_1\in\Gamma(TM_1)$ and 
$X_2\in \Gamma(TM_2)$. We will prove that $X_1 \in \aut(M_1)$ and $X_2\in \aut(M_2)$. First note that
$X_1$ and $X_2$ are Killing vector fields.  

Recall that
\[\mathcal{L}_X J_{ij} =\sum_{k<l} \alpha^{(ij)}_{kl}\, J_{kl},\]
and the endomorphisms \cite{Moroianu-Semmelmann}
\[J_{12}^\pm=\pm {1\over 2} (J_{14}\pm J_{23}),
\hspace{1 cm} J_{31}^\pm=\pm {1\over 2} (J_{13}\mp J_{24}),
\hspace{1 cm} J_{23}^\pm=\pm {1\over 2} (J_{12}\pm J_{34}),\]
where $J_{kl}^-$ and $J_{kl}^+$ vanish on  $M_1$ and $M_2$ respectively.
Let $Z=Z_1+Z_2$ with $Z_1\in\Gamma(TM_1)$
and 
$Z_2\in \Gamma(TM_2)$,
\begin{eqnarray*}
(\Lie_{X_1} J_{ij}^+)(Z_1)_p &=&\Lie_{X_1}(J_{ij}^+(Z_1))_p-J_{ij}^+(\Lie_{X_1}Z_1)_p\in T_pM_1,\\
(\Lie_{X_1} J_{ij}^+)(Z_2)_p &=&\Lie_{X_1}(J_{ij}^+(Z_2))_p-J_{ij}^+(\Lie_{X_1}Z_2)_p=0,\\
(\Lie_{X_2} J_{ij}^+)(Z_1)_p &=&\Lie_{X_2}(J_{ij}^+(Z_1))_p-J_{ij}^+(\Lie_{X_2}Z_1)_p=0,\\
(\Lie_{X_2} J_{ij}^+)(Z_2)_p &=&\Lie_{X_2}(J_{ij}^+(Z_2))_p-J_{ij}^+(\Lie_{X_2}Z_2)_p=0. 
\end{eqnarray*}
i.e. $\Lie_{X} J_{ij}^+=\Lie_{X_1} J_{ij}^+\in \End(TM_1)$. Similarly,
$\Lie_{X} J_{ij}^-=\Lie_{X_2} J_{ij}^-\in \End(TM_2)$. 
Now consider, for instance,
\begin{eqnarray*}
\Lie_{X} J_{12}^+
  &=&{1\over 2}\Lie_{X}( J_{14}+J_{23})\\
  &=&{1\over 2}\sum_{k<l} \alpha^{(14)}_{kl}\, J_{kl}+\sum_{k<l} \alpha^{(23)}_{kl}\, J_{kl}\\
  &=&{1\over 2}\sum\beta_{st}J_{st}^++\sum\gamma_{st}J_{st}^-, 
\end{eqnarray*}
for some functions $\beta_{st}$ and $\gamma_{st}$. Since $\Lie_{X} J_{12}^+\in \End(TM_1)$, all the
coefficients $\gamma_{st}=0$. Therefore 
\[\Lie_{X_1} J_{12}^+=\sum\beta_{st}J_{st}^+,\]
By similar calculations $\Lie_{X_1}
J_{ij}^+=\sum\beta^{(ij)}_{st}J_{st}^+$ and $\Lie_{X_2} J_{ij}^-=\sum\gamma^{(ij)}_{st}J_{st}^-$, i.e.
$X_1\in \aut(M_1)$  and $X_2\in \aut(M_2)$.
Now let $m_i=\dim(M_i)/4$, $i=1,2$.  Since
\begin{eqnarray*}
\dim(\Aut(M))&=&{{2m_1+1} \choose 2} +{{2m_2+1} \choose 2}+6+4m_1+4m_2  \\
   &=& \dim(\Aut(M_1))+\dim(\Aut(M_2)),\\
\dim(\Aut(M_i))&\leq& {{2m_i+1} \choose 2}+3+4m_i
\end{eqnarray*}
we must have  
\[\dim(\Aut(M_i))={{2m_i+1} \choose 2}+3+4m_i.\]
Therefore, the dimensions of the automorphism groups of $M_1$ and $M_2$ are maximal.

For $r\geq 5$, the symmetric spaces carrying a parallel even
Clifford structure were classified in \cite{Moroianu-Semmelmann} and are listed in Table 4. 
We claim that the dimension of the automorphism group of each space listed is maximal
Indeed, let
$M=G/K$ be one of the spaces in Table 4, where $G$ is the group of isometries of $M$. 
We need to prove that every Killing vector is an infinitesimal
automorphism of $M$. If $X\in\mathfrak{g}$ is a Killing vector field,  
\begin{eqnarray*}
(\Lie_X J_{ij})(Z)&=&\Lie_{X}(J_{ij}(Z))-J_{ij}(\Lie_{X}Z)\\
&=&\nabla_X(J_{ij}(Z))-\nabla_{J_{ij}(Z)}X-J_{ij}(\nabla_XZ)+J_{ij}(\nabla_ZX)\\
&=&(\nabla_XJ_{ij})(Z)-[\nabla X,J_{ij}](Z)\\
&=&\left(\sum_{k<l} a_{kl}^{(ij)} J_{kl}\right)(Z)-[\nabla X,J_{ij}](Z)
\end{eqnarray*}
since the almost even-Clifford hermitian structure is parallel.

Recall that
\[\mathfrak{g}=\mathfrak{k}+\mathfrak{m},\]
and at a point $p\in M$, 
\begin{eqnarray*}
\mathfrak{m}&\cong & T_pM,
\end{eqnarray*}
and from the list of possible spaces
\begin{eqnarray*}
\mathfrak{k} &\cong& \mathfrak{a}_p \cong C_{\mathfrak{so}(N)}(\widehat{\mathfrak{spin}({r})})\oplus
\widehat{\mathfrak{spin}({r})}. 
\end{eqnarray*}
On the other hand, since $M$ is symmetric,
\[\mathfrak{g} \cong \mathfrak{b}_p\oplus T_pM,\]
so that 
$\mathfrak{b}_p\cong \mathfrak{k}\cong C_{\mathfrak{so}(N)}(\widehat{\mathfrak{spin}({r})})\oplus
\widehat{\mathfrak{spin}({r})}$ and 
\[[\nabla X,J_{ij}]=\sum_{k<l}
b_{kl}^{(ij)}J_{kl},\] 
i.e.
\[\mathcal{L}_X J_{ij}=\sum_{k<l} (a_{kl}^{(ij)}+b_{kl}^{(ij)})J_{kl}.\]
  \qd

\begin{theo}\label{theo:gap theorem} Let $M$ be a $N$-dimensional rank $r\geq3$, almost even-Clifford
hermitian
manifold and $p\in M$. Assume the following constraints:
\[  
\begin{array}{|c|c|c|c|}
 \hline
  r \mbox{\ {\rm (mod 8)} }  &N &\mbox{ constraint}&\mbox{ extra constraint}  
\tstrut\\
\hline
 0 & d_r(m_1+m_2)& m_1\geq m_2 > {r\choose 2}+1& m_1\equiv m_2 \equiv 0 \mbox{\ {\rm (mod 2)} }\tstrut\\
 & & \mbox{\rm or } m_2\geq m_1> {r\choose 2}+1&\\
 \hline
 1,7 & d_rm& m>{r\choose 2}+1& m \equiv 0 \mbox{\ {\rm (mod 2)} }\tstrut\\
\hline
 2,6 & d_rm& m>{1\over 2}{r\choose 2}+{1\over 2}& m \equiv 0 \mbox{\ {\rm (mod 2)} }\tstrut\\
\hline
 3,5 & d_rm& m> {1\over 4}{r\choose 2} +1&\tstrut\\
\hline
 4 & d_r(m_1+m_2)& m_1\geq m_2 > {1\over 4}{r\choose 2}+1&\tstrut\\
 & & \mbox{\rm or } m_2\geq m_1> {1\over 4}{r\choose 2}+1&\tstrut\\
\hline
  \end{array}
\]
\centerline{\rm Table 5}
If $\dim(\Aut(M))$ is not maximal, then
\[\dim(\Aut(M))<d_{max}-{r\choose2}.\]
\end{theo}
{\em Proof}. 
Suppose that  
\[d_{max}-{r\choose2}\leq \dim(\Aut(M))<d_{max}.\]
At a point $p\in M$, the isotropy group satisfies
\begin{eqnarray*}
\dim(A_p)
  &\geq& \dim(\Aut(M))-N \\
  &\geq& d_{max}-{r\choose 2}-N \\
  &\geq& d_C:=\dim(C_{{SO}(N)}({{Spin}({r})})).
\end{eqnarray*}
The Lie algebra $\mathfrak{a}_p$ of $A_p$ maps one-to-one into
$C_{\mathfrak{so}(N)}(\widehat{\mathfrak{spin}({r})})\oplus \widehat{\mathfrak{spin}({r})}$ since a Killing
vector field $X$ is determined by its values $X_p$ and $(\nabla X)_p$. Consider the
compositions of this map with the projections to the two factors
\[\rho_1:\mathfrak{a}_p\lra
C_{\mathfrak{so}(N)}(\widehat{\mathfrak{spin}({r})}),\quad\quad\rho_2:\mathfrak{a}_p\lra
\widehat{\mathfrak{spin}(r)}.\]
The subalgebra $\rho_1(\mathfrak{a}_p)$ is either equal to
$C_{\mathfrak{so}(N)}(\widehat{\mathfrak{spin}({r})})$
or is contained in a proper maximal subalgebra of $C_{\mathfrak{so}(N)}(\widehat{\mathfrak{spin}({r})})$.

If $r\not \equiv \pm 2 \mbox{\ {\rm (mod 8)}}$, the maximal dimension $d_M$ of a proper maximal subalgebra
of
$C_{\mathfrak{so}(N)}(\widehat{\mathfrak{spin}({r})})$ is
given in the following table (see \cite{mann}), 
\[  
\begin{array}{|c|c|c|c|}
 \hline
  r \mbox{\ {\rm (mod 8)} }  &N &d_M &d_C
\tstrut\\
 \hline
  0 & d_r(m_1+m_2)& \max\left\{{m_1-1\choose 2}+{m_2\choose 2},{m_1\choose 2}+{m_2-1\choose 2}\right\}
&{m_1\choose 2}+{m_2\choose 2}\tstrut\\
 \hline
 1,7 & d_rm&{m-1\choose 2}&{m\choose 2}\tstrut\\
\hline
 3,5 & d_rm& {2m-1\choose 2}+3& {2m+1\choose 2}\tstrut\\
\hline
 4 & d_r(m_1+m_2)& \max\left\{{2m_1-1\choose 2}+3+{2m_2+1\choose 2},{2m_1+1\choose 2}+{2m_2-1\choose
2}+3\right\}
&{2m_1+1\choose 2}+{2m_2+1\choose 2}\tstrut\\
\hline
  \end{array}
\]
\centerline{\rm Table 6}
Thus,  due to the constraints on the multiplicities $m,m_1,m_2$, if $\rho_1(\mathfrak{a}_p)$ is contained in a
proper subalgebra of $C_{\mathfrak{so}(N)}(\widehat{\mathfrak{spin}({r})})$,
\begin{eqnarray*}
  d_C
   &>& 
  d_M + {r\choose 2}\\
   &\geq&
  \dim(\rho_1(\mathfrak{a}_p)) +  \dim(\rho_2(\mathfrak{a}_p)) \\
   &\geq&
  \dim(\mathfrak{a}_p) \\
   &\geq&
   d_C,
\end{eqnarray*}
which is a contradiction. Thus
\[\rho_1(\mathfrak{a}_p) =C_{\mathfrak{so}(N)}(\widehat{\mathfrak{spin}({r})}), \]
and
\[\mathfrak{a}_p\cong
C_{\mathfrak{so}(N)}(\widehat{\mathfrak{spin}({r})})\oplus \mathfrak{K}\subset\mathfrak{so}(N),\]
where
\[\mathfrak{K}:=\ker(\rho_1|_{\mathfrak{a}_p})\subset \ker(\rho_1)= \mathfrak{spin}(r).\]
Therefore  $A_p=C_{SO(N)}(Spin(r))\cdot K$, where $K$ is a closed subgroup of $Spin(r)$. The extra assumptions
on $m$, $m_1$ and $m_2$
imply $-1\in C_{SO(N)}(Spin(r))$ and $-1\in A_p$. 
Thus, there is an element $g\in A_p$ whose derivative
$dg_p=-{\rm Id}_{T_pM}$ in the isotropy representation of $A_p$ on $T_pM$. 
In other words, the automorphism
$g$ is a (global) symmetry at $p$ and $M$ is symmetric. 
Since these symmetries generate the translations along geodesics, $M$ has a
transitive group of automorphisms, not just isometries. As in the proof of Proposition
\ref{prop:sym-parallel}, this implies the almost even-Clifford hermitian structure is parallel. By
arguments similar to those in the proof of Theorem \ref{main-theorem}, we have $\dim(\Aut(M))=d_{max}$, which
is again a contradiction.

For $r \equiv \pm2 \mbox{\ {\rm (mod 8)}}$, $\dim(A_p)\geq d_C$ will happen if
$\mathfrak{su}(m)\subset \rho_1(\mathfrak{a}_p)$ . Hence $A_p=H\cdot K$, where $H$ is some subgroup of
$C_{SO(N)}(Spin(r))$ containing
$SU(m)$ and $K$ is some closed subgroup of $Spin(r)$. 
Since we are assuming $m$ is even, $-1\in SU(m)$ and $-1\in A_p$. Therefore $M$ is symmetric and again,
the proofs of Proposition \ref{prop:sym-parallel} and Theorem \ref{main-theorem} imply
$\dim(\Aut(M))=d_{max}$.
\qd

{\bf Remark}. The constraints in the previous theorem are given  in order to ensure that $-1\in
A_p$. If we were to relax them or change them, a more detailed analysis of the possible subgroups $K\in
Spin(r)$ would be needed.

{\small
\renewcommand{\baselinestretch}{0.5}
\newcommand{\bi}{\vspace{-.05in}\bibitem} }

\end{document}